# Solvable Polynomials Over The Gaussian Field Q(*i*)

Nicholas Phat Nguyen[1]

**Abstract.** In this paper, I outline how the results from my previous paper (*A Congruence Property of Solvable Polynomials*) could be generalized with some modification to polynomials with coefficients in the Gaussian field **Q**(*i*).

**A. IMMAGINARY QUADRATIC FIELDS OF CLASS NUMBER ONE.** Let *L* be an imaginary quadratic number field and $\mathcal{R}$ its ring of integers. If $\mathcal{R}$ is a principal ideal domain, then there are some similarities between such a quadratic field *L* and the field **Q** of rational numbers.

- Each prime ideal of $\mathcal{R}$ is principal. Numbers of $\mathcal{R}$ that generate prime ideals are also called prime numbers.
- Class field theory tells us that for each integer $\alpha$ in the ring $\mathcal{R}$, there is a Galois extension $L_\alpha$ of *L* such that the Galois group Gal($L_\alpha$ /*L*) of the extension is isomorphic to the group $(\mathcal{R}/\alpha\mathcal{R})^*/\mathcal{R}^*$ of invertible residue classes in the quotient ring $\mathcal{R}/\alpha\mathcal{R}$ modulo the residue classes represented by the units (the roots of unity) in the ring $\mathcal{R}$. The extension $L_\alpha$ is known as the ray class field modulo $\alpha$ of *L*, and it is analogous to the cyclotomic extension **Q**($\zeta_n$) of **Q**.
- Any prime ideal in *L* that ramifies in $L_\alpha$ must divide $\alpha$, and each abelian extension of *L* is contained in one of the ray class fields $L_\alpha$.

There are exactly 9 imaginary quadratic fields whose ring of integers is a principal ideal domain.[2] These are known as imaginary quadratic fields of class number one (because the ring of integers of such a field is a principal ideal domain if and only if the ideal class group of that field is trivial).

---

[1] E-mail address: nicholas.pn@gmail.com

[2] They are the imaginary quadratic fields with discriminants -3, -4, -7, -8, -11, -19, -43, -67, -163.

Given such an imaginary quadratic field $L$ of class number one, the ring $\mathcal{R}$ is a full lattice in the complex plane $\mathbf{C}$, and quotient group $\mathbf{C}/\mathcal{R}$ is isomorphic (as a complex Lie group) to an elliptic curve $\mathcal{E}$ defined over the field $\mathbf{Q}$ of rational numbers.

The ring of endomorphisms of the elliptic curve $\mathcal{E}$ is isomorphic to $\mathcal{R}$, so that $\mathcal{E}$ is said to have complex multiplication. For each integer $\alpha$ in the ring $\mathcal{R}$, the $\alpha$-torsion points of the $\mathcal{E}$ (regarded as an $\mathcal{R}$-module) form a free $\mathcal{R}/\alpha\mathcal{R}$ module of rank one, and its automorphism group is naturally isomorphic to the group $(\mathcal{R}/\alpha\mathcal{R})^*$ of invertible residue classes in the quotient ring $\mathcal{R}/\alpha\mathcal{R}$, i.e., the group of units of $\mathcal{R}/\alpha\mathcal{R}$. If we denote $N(\alpha)$ for the norm of $\alpha$ in the quadratic field $L$, then there are $N(\alpha)$ $\alpha$-torsion points in $\mathcal{E}$.

Now take a Weber function $\omega$ defined on $\mathcal{E}$. By definition $\omega$ is a rational function on $\mathcal{E}$ with coefficients in $\mathbf{Q}$, and $\omega$ is invariant under automorphisms of $\mathcal{E}$ over $\mathbf{Q}$. The theory of complex multiplication tells us that the ray class field $L_\alpha$ can be obtained by adjoining to $L$ the values of the Weber function $\omega$ evaluated at the $\alpha$-torsion points of $\mathcal{E}$.

By multiplying $\omega$ with a rational integer as needed, we can assume that the values $\omega(T)$, where $T$ runs through all the $\alpha$-torsion points of $\mathcal{E}$ (for a specific $\alpha$ under consideration), are algebraic integers. For an integer $u$ in $\mathcal{R}$ prime to $\alpha$, if $(u, L_\alpha/L)$ is the element of $\text{Gal}(L_\alpha/L)$ defined by $u$ under the Artin reciprocity map, then the Galois action of $(u, L_\alpha/L)$ can be described as follows, assuming $u$ is not divisible by any prime $p$ over which the elliptic curve $\mathcal{E}$ has bad reduction modulo $p$. If $T$ is an $\alpha$-torsion point of $\mathcal{E}$, then the Galois action by $(u, L_\alpha/L)$ transforms $\omega(T)$ to $\omega(uT)$, where $uT$ represents the transform of $T$ under the complex multiplication by $u$ on $\mathcal{E}$. Note that for a prime number $\pi$ of $L$ not dividing $\alpha$, the action of $(\pi, L_\alpha/L)$ on the residue class $\omega(T)$ modulo any prime ideal of $L$ above $\pi$ is just the Frobenius automorphism with $\mathcal{R}/\pi\mathcal{R}$ as the fixed field.

Hence, the ray class field adjunction for $L$ through torsion points of an elliptic curve model for $\mathbf{C}/\mathcal{R}$ and the nice Galois action in this case are similar to the ray class field adjunction of **Q** through torsion points of the unit circle and the Galois action on roots of unity. Given these similarities it is tempting to consider whether we can generalize the results in my previous paper (*A Congruence Property of Solvable Polynomials*) to solvable polynomials over any imaginary quadratic field of class number one. However, there are some technical hurdles, the biggest of which is whether we can find a suitable Weber function $\omega$ for which all the values $\omega(T)$ are algebraic integers, where $T$ runs over torsion points relative to <u>all</u> integers $\alpha$ in $\mathcal{R}$ that are not divisible by the finite number of primes where $\mathcal{E}$ has bad reduction.

In case $L = \mathbf{Q}(i)$, it happens that we can find a rational function on $\mathcal{E}$ that, while not being a Weber function, meet our requirements above and play the same role as the exponential function for cyclotomic extensions over **Q**.

**B. THE GAUSSIAN FIELD Q($i$) AND THE LEMNISCATE FUNCTION.** For the field $L = \mathbf{Q}(i)$ and $\mathcal{R} = \mathbf{Z}[i]$, we can use the lemniscate function. For real arguments, this function is the inverse of the arc length integral on the Bernoulli lemniscate given by the equation $r^2 = \cos(2\theta)$ in polar coordinates. That function can be extended to an elliptic function $\varphi(z)$ with period lattice $\mathcal{L} = \mathbf{Z}(1+i)\varpi + \mathbf{Z}(1-i)\varpi$, where $\varpi$ is one-quarter of the length of the lemniscate, i.e., $\varpi = \int dt/(1-t^4)^{1/2}$, the integral being taken from $t=0$ to $t=1$. The period lattice $\mathcal{L}$ is proportional to $\mathbf{Z}[i]$ and therefore defines a **C**-elliptic curve isomorphic to $\mathbf{C}/\mathcal{R}$.

We use the curve $Y^2 = 4X^3 + X$ as a model for $\mathbf{C}/\mathcal{L}$. That curve is the elliptic curve parametrized by the Weierstrass elliptic functions associated with the period lattice $\mathbf{C}/\mathcal{L}$, and so the lemniscate function can also be viewed as a rational function on that elliptic curve.[3] See [5] at Lemma 2.2. The discriminant of the curve is divisible by only the rational prime 2, and so this curve has good reduction at all odd prime numbers.

---

[3] The lemniscate function is not a Weber function because it is not invariant under the automorphisms of that elliptic curve.

The field $\mathbf{Q}(i)$ has discriminant $-4$, and 2 is the only rational prime number that ramifies in $\mathbf{Q}(i)$. We have $2 = (1 + i)(1 - i)$ with $(1 + i) = i(1 - i)$ being a prime number in $\mathbf{Q}(i)$. The group $\mathcal{R}^*$ of units in $\mathbf{Z}[i]$ consists of the four roots of unity $\{1, -1, i, -i\}$.

We will call integers in $\mathbf{Z}[i]$ that are not divisible by the prime number $(1 + i)$ the odd integers. Note that the model curve $Y^2 = 4X^3 + X$ of $\mathbf{C}/\mathcal{L}$ has good reduction at any odd prime number of $\mathcal{R}$.

For an odd integer $\beta$, the ray class field modulo $\alpha = 2(1 + i)\beta$ has Galois group isomorphic to $(\mathcal{R}/\alpha\mathcal{R})^*/\mathcal{R}^*$. Note that $(\mathcal{R}/\alpha\mathcal{R})^*$ is naturally isomorphic to $(\mathcal{R}/\beta\mathcal{R})^*$ x $(\mathcal{R}/2(1+i)\mathcal{R})^*$, and the four roots of unity in $\mathcal{R}^*$ correspond to the four invertible residue classes in $\mathcal{R}/2(1+i)\mathcal{R})$. Accordingly, $(\mathcal{R}/\alpha\mathcal{R})^*/\mathcal{R}^*$ is naturally isomorphic to $(\mathcal{R}/\beta\mathcal{R})^*$.

In a 2014 paper, David Cox and Trevor Hyde studied the lemniscate function and its Galois property with respect to the ray class fields $L_\alpha$ of $\mathbf{Q}(i)$ modulo $\alpha = 2(1 + i)\beta$, where $\beta$ is an odd integer. See [5]. They found the following amazing facts.

- All values of the complex lemniscate function $\varphi(z)$ at $\beta$-torsion points of $\mathbf{C}/\mathcal{L}$ are distinct algebraic integers.
- Recall that the $\beta$-torsion points of $\mathbf{C}/\mathcal{L}$ form a free $(\mathcal{R}/\beta\mathcal{R})$-module of rank one. Pick any generator $S$ of that module. The extension $L_\alpha$ is generated by the values $\varphi(\lambda S)$ at the points $\lambda S$, with $\lambda$ running through all the invertible elements of $\mathcal{R}/\beta\mathcal{R}$.
- Let $\rho$ be the Galois automorphism of the ray class field $L_\alpha$ corresponding to the residue class of an odd prime $\pi$ in $(\mathcal{R}/\beta\mathcal{R})^*$ under the Artin reciprocity map, then $\rho(\varphi(U)) = \varphi(\pi U)$ for any $\beta$-torsion point U of $\mathbf{C}/\mathcal{L}$.

Accordingly, the ray class field $L_\alpha$ can be obtained by adjoining to $\mathbf{Q}(i)$ any number $\varphi(\lambda S)$, where $\lambda$ is invertible in $\mathcal{R}/\beta\mathcal{R}$. We will call the numbers $\varphi(\lambda S)$ primitive $\beta$-torsion values in analogy with the cyclotomic case. All such primitive $\beta$-torsion values are roots of an irreducible monic polynomial with coefficients in $\mathcal{R}$, which Cox and Hyde call the $\beta^{th}$ lemnatomic polynomial. We will denote that lemnatomic polynomial as $\Lambda_\beta$. Its degree is the same as the degree of $L_\alpha$ over $\mathbf{Q}(i)$.

**C. LEMNATOMIC POLYNOMIALS.** The goal of this section is to show that for any odd prime number $\pi$ in $\mathcal{R} = \mathbf{Z}[i]$ that does not divide the odd integer $\beta$, the $\beta^{\text{th}}$ lemnatomic polynomial $\Lambda_\beta$ is separable modulo $\pi$.

Again let $\alpha = 2(1 + i)\beta$ and denote by $L_\alpha$ the ray class field modulo $\alpha$ of $\mathbf{Q}(i)$. Let $J$ be the ring of integers of $L_\alpha$. From class field theory, we know that the only prime numbers in $\mathbf{Q}(i)$ that ramify in $L_\alpha$ are $(1 + i)$ and the prime divisors of $\beta$. Accordingly, if $\pi$ is an odd prime number that does not divide $\beta$, then the principal ideal $\pi J$ is a product of distinct prime ideals in $J$ with no repetition.

Let $m$ be a prime ideal in $J$ containing $\pi$. We can embed $L_\alpha$ as a dense subring in an unramified extension $W$ of the $\pi$-adic completion of $\mathbf{Q}(i)$, so that $m$ is also a dense subset of the maximal ideal in the ring of integers of $W$. The extension $W$ is generated by a primitive $\beta$-torsion value $\varphi(S)$. The residue class of $\varphi(S)$ modulo $\pi$ generates the residue field of $\pi$ in $W$. Moreover, because $W$ is unramified over $\mathbf{Q}(i)_\pi$, $\varphi(S)$ also generates the ring of integers of $W$, i.e., the ring of integers of $W$ is $\mathcal{R}_\pi[\varphi(S)]$.

**Lemma 1.** We have $\pi J + \mathcal{R}[\varphi(S)] = J$.

*Proof.* Both sides are $\mathcal{R}$ modules, and it is sufficient to show that the localization of each side at every prime ideal of $\mathcal{R}$ is the same.

For any prime ideal of $\mathcal{R}$ that does not contain $\pi$, the localizations of both sides are clearly the same. So what remains is to determine the localizations at $\pi$ for both sides. Both sides are finitely-generated modules over the localization of $\mathcal{R}$ at $\pi$. They are the same if and only if they give us the same modules when we pass to their $\pi$-adic completions. But in light of our preceding discussion, the $\pi$-adic completions of both sides give us the ring of integers in $W$. ∎

**Lemma 2.** *The quotient ring $J/\pi J$ is isomorphic to $(\mathcal{R}/\pi\mathcal{R})[X] / \Lambda_\beta(X)$.*

*Proof.* By Lemma 1, we know $J/\pi J$ is naturally isomorphic to $\mathcal{R}[\varphi(S)] / (\mathcal{R}[\varphi(S)] \cap \pi J)$. Now $\mathcal{R}[\varphi(S)] \cap \pi J \supset \pi \mathcal{R}[\varphi(S)]$, so there is a surjection:

$$\mathcal{R}[\varphi(S)]/\pi\mathcal{R}[\varphi(S)] \twoheadrightarrow J/\pi J.$$

Note that $\mathcal{R}[\varphi(S)]$ is naturally isomorphic to $\mathcal{R}[X]/\Lambda_\beta(X)$, so $\mathcal{R}[\varphi(S)]/\pi\mathcal{R}[\varphi(S)]$ is naturally isomorphic to $(\mathcal{R}/\pi\mathcal{R})[X] / \Lambda_\beta(X)$.

Because $J/\pi J$ is a vector space over $\mathcal{R}/\pi\mathcal{R}$ of dimension $\deg(L_\alpha : \mathbf{Q}(i)) = \deg \Lambda_\beta(X)$, same as the dimension of $(\mathcal{R}/\pi\mathcal{R})[X] / \Lambda_\beta(X)$, the above surjection must be bijective. Therefore the quotient ring $J/\pi J$ is isomorphic to $(\mathcal{R}/\pi\mathcal{R})[X] / \Lambda_\beta(X)$. ∎

**Proposition 1.** *The the $\beta^{th}$ lemnatomic polynomial $\Lambda_\beta(X)$ is separable modulo $\pi$ for any odd prime number $\pi$ in $\mathcal{R}$ that does not divide $\beta$.*

*Proof.* Lemma 2 above says that $(\mathcal{R}/\pi\mathcal{R})[X]/\Lambda_\beta(X)$ is isomorphic to $J/\pi J$, a reduced ring (i.e., having no nilpotent elements). Accordingly, $\Lambda_\beta(X)$ is separable modulo $\pi$.[4] ∎

Note that for any extension $K$ of $\mathbf{Q}(i)$, the above proposition implies that $\Lambda_\beta(X)$ is separable modulo any prime ideal in that extension that does not divide $(1 + i)\beta$.

**D. A CRITERION FOR IRREDUCIBILITY.** By the results of Cox and Hyde, the lemnatomic polynomials are irreducible over the field $\mathbf{Q}(i)$. We can now give a criterion for a lemnatomic polynomial to be irreducible over an extension $K$ of $\mathbf{Q}(i)$, generalizing our propositions about cyclotomic polynomials from [6].

**Proposition 2:** *Let $K$ be an extension of $\mathbf{Q}(i)$. Let $\mathcal{R} = \mathbf{Z}[i]$ be the ring of integers in $\mathbf{Q}(i)$, and let $\beta$ be an odd integer in $\mathcal{R}$. If each residue class in a generating set of the unit group $(\mathcal{R}/\beta\mathcal{R})^*$ contains an odd prime number of $Q(i)$ that is semi-split in $K$, then the $\beta^{th}$ lemnatomic polynomial $\Lambda_\beta(X)$ is irreducible over K.*

---

[4] Over a finite field, or any perfect field, a polynomial is separable if and only if it has no repeated irreducible factors.

*Proof.* The proof of this proposition is similar to the case over $\mathbb{Q}$. However, because of the higher level of abstraction, it is important for us to go over the proof in this new context to make sure that all details work.

Each irreducible polynomial in $K[X]$ corresponds to an orbit in a given separable closure of $K$ under the action of the absolute Galois group of $K$. All elements in that orbit are the roots of the irreducible polynomial. If $\Lambda_\beta(X)$ is not irreducible over $K$, then there must be two different roots of $\Lambda_\beta(X)$, say $u$ and $v$, that belong to two distinct Galois orbits over $K$. The integral elements $u$ and $v$ must then have distinct minimal polynomials $g$ and $h$ over $K$ whose product divides $\Lambda_\beta(X)$ in $\mathcal{R}[X]$.

Let $A$ be the ring of integers of $K$ and let $B$ be the ring of integers in the extension $K(u, v)$. For each maximal ideal $\mathfrak{m}$ of $A$ that does not divide $(1 + i)\beta A$, $\Lambda_\beta(X)$ mod $\mathfrak{m}$ is separable by Proposition 1. So if we consider any maximal ideal $\wp$ of $B$ in the prime ideal factorization of $\mathfrak{m}B$, the images $u$ mod $\wp$ and $v$ mod $\wp$ in the residue field $B/\wp$ must belong to different Galois orbits over $A/\mathfrak{m}$, because they are the roots of polynomials $g$ mod $\mathfrak{m}$ and $h$ mod $\mathfrak{m}$ whose product divides the separable polynomial $\Lambda_\beta(X)$ mod $\mathfrak{m}$.

Because both $u$ and $v$ are roots of $\Lambda_\beta(X)$, we must have $u = \varphi(U)$ and $v = \varphi(V)$, where $\varphi$ is the complex lemniscate function, and $U$ and $V$ are some primitive $\beta$-torsion points in $\mathbb{C}/\mathcal{L}$. Because the group of $\beta$-torsion points is a free $\mathcal{R}/\beta\mathcal{R}$ module of rank one, there is an element $t$ in $(\mathcal{R}/\beta\mathcal{R})^*$ such that $V = tU$. For our purpose, we can assume that $t$ belongs to the given generating set for the group $(\mathcal{R}/\beta\mathcal{R})^*$. Indeed, if both $U$ and $tU$ give us two values $\varphi(U)$ and $\varphi(tU)$ in the same Galois orbit over $K$ regardless of which $U$ we choose and which $t$ in the generating set we choose, then all the values $\varphi(U)$ as $U$ runs through the $\beta$-torsion points must be in the same Galois orbit over $K$, contrary to our initial assumption.

By hypothesis, we can find an odd prime number $\pi$ in the ring $\mathcal{R}$ that is semi-split in $K$ and congruent to $t$ mod $\beta\mathcal{R}$. This odd prime $\pi$ does not divide $\beta$. So for each maximal ideal $\mathfrak{m}$ in $A$ that contains $\pi$, $\Lambda_\beta(X)$ is separable over $A/\mathfrak{m}$. Consequently, for any maximal ideal $\wp$ of $B$ in the factorization of $\mathfrak{m}B$, the elements $u$ mod $\wp$ and $v$ mod $\wp$ will be in distinct Galois orbits over $A/\mathfrak{m}$.

If $\rho$ is the Galois automorphism of the ray class field $L_\alpha$ corresponding to the residue class of an odd prime $\pi$ in $(\mathcal{R}/\beta\mathcal{R})^*$ under the Artin reciprocity map, the Galois theory of the lemniscate by Cox and Hyde tells us that $\rho(\varphi(U)) = \varphi(\pi U) = \varphi(V)$. Moreover, by class field theory, $\rho$ is the Frobenius automorphism associated with $\pi$. That means for the integers $u = \varphi(U)$ and $v = \varphi(V)$, $v \bmod \wp$ is the transform of $u \bmod \wp$ under the map which raises each element in the residue field $B/\wp$ to the exponent $N(\pi)$, where the norm $N(\pi)$ of $\pi$ in $\mathbf{Q}(i)$ represents the number of elements in the residue field $\mathcal{R}/\pi\mathcal{R}$.

Because by hypothesis $\pi$ is semi-split in $K$, we can choose a maximal ideal $m$ of $A$ containing $\pi$ such that $A/m$ is isomorphic to $\mathcal{R}/\pi\mathcal{R}$. That means $u \bmod \wp$ and $v \bmod \wp$ must belong to the same Galois orbit over $A/m$ because the Frobenius automorphism in question is a generator for any finite Galois extension of $\mathcal{R}/\pi\mathcal{R}$.

However, we saw earlier that $u \bmod \wp$ and $v \bmod \wp$ must be in different Galois orbits over $A/m$ because $u$ and $v$ are assumed to be in different Galois orbits over $K$. This contradiction shows that all the roots of $\Lambda_\beta(X)$ must be in one Galois orbit over $K$. Therefore $\Lambda_\beta(X)$, must be irreducible over $K$. ∎

For any odd integer $\delta$ relatively prime to $\beta$, the ray class field modulo $2(1 + i)\delta$ would meet the above condition. The odd prime numbers $\pi$ in $\mathcal{R}$ that are congruent to 1 modulo $2(1 + i)\delta\mathcal{R}$ are semi-split in that ray class field. The Chinese remainder and Chebotarev density theorems tell us that there are always odd prime numbers $\pi$ in $\mathcal{R}$ (indeed infinitely many) that satisfy both the congruence $p \equiv 1 \pmod{2(1 + i)\delta\mathcal{R}}$ and the congruence $p \equiv t \pmod{\beta\mathcal{R}}$ for each invertible residue class $t \bmod \beta\mathcal{R}$. These prime numbers would be semi-split in the given ray class field and would also belong to the residue class $t \bmod \beta\mathcal{R}$. Accordingly, the $\beta^{th}$ lemnatomic polynomial $\Lambda_\beta(X)$ is irreducible over the ray class field modulo $2(1 + i)\delta$ of $\mathbf{Q}(i)$.

**Proposition 3**. *Let $K$ be an extension of $\mathbf{Q}(i)$ with odd discriminant relative to $\mathbf{Q}(i)$. For $K$ to contain a nontrivial solvable extension of $\mathbf{Q}(i)$, it is necessary that we can find an odd integer $\beta$ in $\mathcal{R} = \mathbf{Z}[i]$ where every generating set of $(\mathcal{R}/\beta\mathcal{R})^*$ has at least a residue class containing <u>no</u> odd prime number of $\mathbf{Q}(i)$ that is semi-split in $K$.*

*Proof.* By the definition of solvable extensions, $K$ contains a nontrivial solvable extension of $\mathbf{Q}(i)$ if and only if $K$ contains a nontrivial abelian extension of $\mathbf{Q}(i)$. From class field theory, we know any abelian extension of $\mathbf{Q}(i)$ is contained in some ray class field extension. So $K$ contains a nontrivial abelian extension of $\mathbf{Q}(i)$ if and only if $K \cap L_\alpha \neq \mathbf{Q}(i)$ for some ray class field $L_\alpha$ modulo an integer α. Moreover, because the discriminant of $K$ relative to $\mathbf{Q}(i)$ is odd, we can take $\alpha$ to be of the form $2(1 + i)\beta$ for some odd integer $\beta$ in $\mathbf{Q}(i)$.

To say that $K \cap L_\alpha \neq \mathbf{Q}(i)$ is equivalent to saying that $\Lambda_\beta(X)$ is reducible over $K$. For $\Lambda_\beta(X)$ to be reducible over $K$, it is necessary that the number field $K$ does not meet the condition of Proposition 2 above, which gives us a sufficient condition for $\Lambda_\beta(X)$ to be irreducible over $K$. In particular, if a number field $K$ fails the condition of Proposition 2 for $\beta$, every generating set of the group $(\mathcal{R}/\beta\mathcal{R})^*$ must include at least a residue class containing no odd prime number semi-split in $K$. ∎

### E. SOLVABLE POLYNOMIALS OVER Q($i$).

Putting things together, we have the following theorem which generalizes our result in [6].

**Theorem :** *Let $\mathcal{R} = \mathbf{Z}[i]$ be the ring of integers in $\mathbf{Q}(i)$. Let $h(X) \in \mathcal{R}[X]$ be a monic polynomial with integer coefficients and whose discriminant is odd, and let $P$ be the set of odd prime numbers $\pi$ in $\mathbf{Q}(i)$ such that $h(X)$ splits completely into distinct linear factors modulo $\pi$. If the Galois group of $h(X)$ over $\mathbf{Q}(i)$ has a nontrivial solvable quotient group, then there is an odd integer $\beta$ in $\mathbf{Q}(i)$ having a nontrivial common factor with the discriminant of $h(X)$, such that every generating set of the unit group $(\mathcal{R}/\beta\mathcal{R})^*$ has at least a residue class containing <u>no</u> prime number $\pi$ in $P$.*

*Proof.* The proof this theorem is based on Propositions 2 and 3 above, and follows essentially the same reasoning as in the case over $\mathbf{Q}$. See [6].

In brief, if the splitting extension $K$ of the polynomial $h$ contains a nontrivial solvable extension of $\mathbf{Q}(i)$, then there is an integer $\alpha$ such that $K \cap L_\alpha \neq \mathbf{Q}(i)$, where $L_\alpha$ is the ray class field modulo $\alpha$ over $\mathbf{Q}(i)$. Because any prime number of $\mathbf{Q}(i)$ ramifies in $K$ if and only if it divides the discriminant of $h$, the discriminant of $K$ must have the same prime factors as the discriminant of $h$, and therefore must also be odd.

We can take $\alpha$ to be of the form $2(1 + i)\beta$ for some odd integer $\beta$ in $\mathbf{Q}(i)$. By Proposition 3 above, every generating set of the unit group $(\mathcal{R}/\beta\mathcal{R})^*$ has at least a residue class containing no odd prime number of $\mathbf{Q}(i)$ that is semi-split in $K$. Based on the same arguments as outlined in [6] and which we do not repeat here, such a residue class would contain no odd prime number $\pi$ of $\mathbf{Q}(i)$ such that $h(X)$ splits completely into distinct linear factors modulo π.

The odd prime numbers of $\mathbf{Q}(i)$ that ramify in $K \cap L_\alpha$ must divide $\beta$ because only those odd prime numbers ramify in the ray class field $L_\alpha$. At the same time, any odd prime number of $\mathbf{Q}(i)$ that ramifies in $K \cap L_\alpha$ must ramify in at least one simple extension obtained by adjoining a root of $h(X)$ to $\mathbf{Q}(i)$, which means such a prime must divide the discriminant of $h(X)$ as well. Accordingly, any odd prime number $\pi$ of $\mathbf{Q}(i)$ that ramifies in $K \cap L_\alpha$ must divide $\beta$ as well as the discriminant of $h(X)$.

Because $\mathbf{Q}(i)$ has class number one, it has no unramified abelian extension other than $\mathbf{Q}(i)$ itself. Because $K \cap L_\alpha$ is a nontrivial abelian extension of $\mathbf{Q}(i)$ and because the discriminant of $K$ is odd, there is at least an odd prime number $\pi$ in $\mathbf{Q}(i)$ that ramifies in the extension $K \cap L_\alpha$, and therefore the number $\beta$ has a nontrivial common factor with the discriminant of $h(X)$. ∎